\documentclass[12pt, english]{article}
\makeatletter
\providecommand{\LyX}{L\kern-.1667em\lower.25em\hbox{Y}\kern-.125emX\@}
\makeatother

\usepackage{Loc}
\usepackage{graphics}
\begin{document}


\begin{center}
\linespread{1}
\Large{\textbf{Localization on Snowflake Domains}}\\
 \vspace{1in}
\Large{\textbf{Britta Daudert}}\\
Department of Mathematics\\
University of California\\
Riverside, CA 92521--0135\\
e-mail address: \textit{britta@math.ucr.edu}\\
\vspace{.5in}
\Large{\textbf{Michel L. Lapidus}}\\
Department of Mathematics\\
University of California\\
Riverside, CA 92521--0135\\
e-mail address: \textit{lapidus@math.ucr.edu}\\
\end{center}

\vspace{.7in}

\begin{abstract}
\noindent The geometric features of the square and triadic Koch
snowflake drums are compared using a position entropy defined on
the grid points of the discretizations (pre-fractals) of the two
domains. Weighted graphs using the geometric quantities are
created and random walks on the two pre-fractals are performed.
The aim is to understand if the existence of narrow channels in
the domain may cause the `localization' of eigenfunctions.

\end{abstract}

\section{Introduction}
\linespread{1.6}

\noindent The term fractal is defined very loosely \cite{F}. Here
are just some ways of characterizing a fractal set  F:

\begin{itemize}
\item F has a fine structure, i.e., detail on an arbitrary small
scale.

\item F is too irregular to be described in traditional geometric
language, both locally and globally.

\item Often, F has some form of self-similarity (maybe strict,
approximate or statistical).


\item Usually, the fractal dimension of F (defined in some way) is
greater than its topological dimension.

\item In most cases, F is defined in a very simple way (e.g.,
recursively).
\end{itemize}

\noindent In short, fractals are objects with irregular geometry.
Such objects can be found everywhere in nature, the most obvious
example being that of a tree structure. The {\bf vibrational}
properties of fractals are of great interest. Questions that arise
naturally are, for example:

\begin{itemize}
\item Why are waves damped much more by fractal rather than by
smooth coastlines?

\item How can we explain the vibrational properties of glass?

\item Why does fractal foliage of trees provide such strong
resistance to the wind?
\end{itemize}

\noindent All these questions are largely unanswered. Fractal
geometry is a means to describing strongly irregular objects and
questions like the ones mentioned above have played and continue
to play an important role in the development of the subject. The
main idea is to use strongly irregular but {\bf deterministic}
objects as a good approximation to the fractal. If the physical
properties of the object of study are related to the character of
the geometry, these deterministic approximations prove to be a
good source of information.

\noindent\subsection{Fractal Drums}

\noindent A fractal drum is the simplest example of a surface
fractal resonator. It is a flat surface bounded by a fractal
frontier.

\vspace{.3in}

\begin{figure}[h]
\label{fig:koch} \caption{An approximation to the triadic Koch
Snowflake domain}\vspace{.7in}\ \ \ \ \ \  \ \ \ \ \ \ \ \ \ \ \ \
\ \ \
 \scalebox{1.2} {\includegraphics{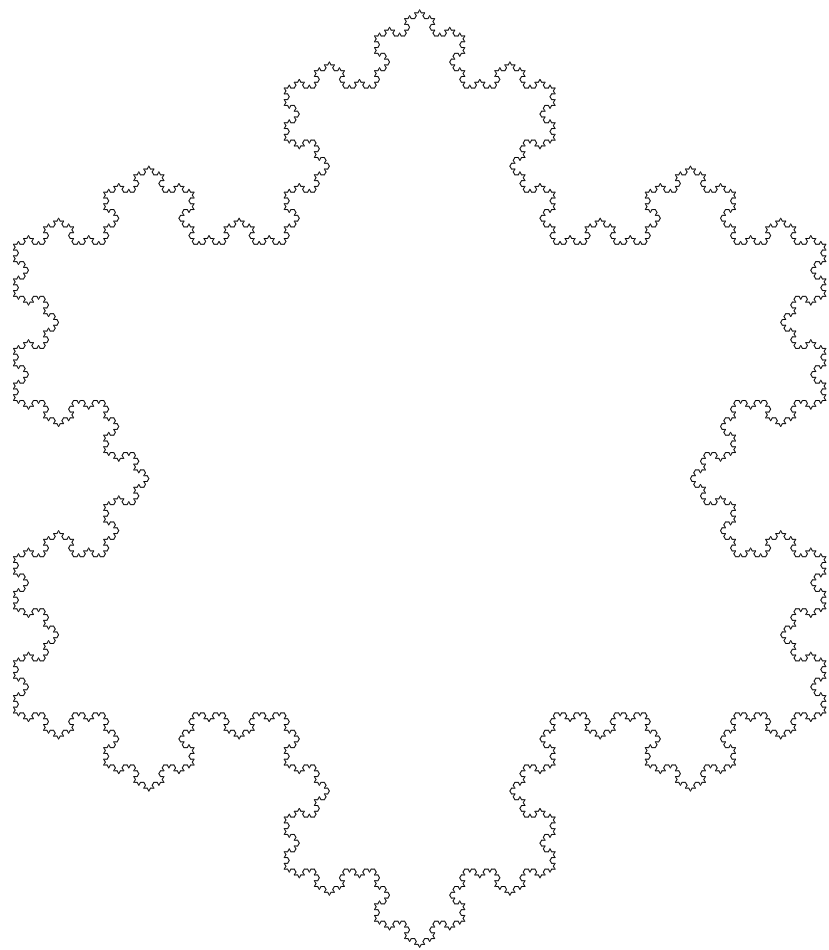}} \label{fig:koch}
\end{figure}

\newpage
\subsection{Vibrational properties of fractal drums}

\noindent The vibrational properties of these structures are of
great theoretical and practical interest \cite{Be1}, \cite{Be2},
\cite{BroCa}, \cite{L1}, \cite{L2}, \cite{L3}, \cite{LNRG},
\cite{LP}, \cite{Sa2}, \cite{EveRRPS}, \cite{SG}, \cite{SGM}.

\vspace{.2in} \noindent An extensive amount of experimental work
on the  study of waves or harmonic oscillations carried by
fractals has been done by the physicist B. Sapoval and his
collaborators \cite{SG}, \cite{SGM}. The {\bf fractal drum} in
their experiments consisted of a stainless steel sheet with the
boundary etched in the pattern of a square Snowflake. A soap
bubble deposited on this object was excited to low frequency
vibrations by a loudspeaker situated above the fractal drum. Their
results showed that the wave motion is strongly damped by the
fractal boundary and that localization may occur. Sapoval
attributed this localization to the existence of narrow channels
in the geometry of the square domain. This phenomenon was not
observed
on the triadic Snowflake domain. \\
In \textit{Snowflake harmonics and computer graphics: numerical
computation of spectra on fractal drums} \cite{LNRG}, the authors
numerically estimated the first 50 eigenvalues and eigenvectors of
the triadic Snowflake domain and tested graphically the results
and conjectures concerning the boundary behavior of the gradient
of
the eigenfunctions.\\
Some mathematical work has been done to confirm Sapoval's results.
Recently, M. L. Lapidus and M. Pang \cite{LP} have studied the
boundary behavior of the Dirichlet Laplacian eigenfunctions and
their gradients on a class of planar domains with fractal
boundary, including the (triangular and square) Koch Snowflake
domains, as well as their polygonal approximations. One of their
main results, specialized to the Koch Snowflake domain, states
that the magnitude of the gradient of the first eigenfunction (or
`first harmonic') `blows up' at infinitely many boundary points,
i.e., the membrane of the Koch snowflake drum exhibits `infinite
stress' near such points. Physically, this corresponds to a
strong damping phenomenon.\\
In this paper, we develop a new approach to investigate
localization phenomena on fractal drums: the discretizations of
the square and triadic Snowflake are considered as weighted
graphs. The weight is assigned to an edge in such a way that it
reflects certain geometric properties of the two nodes connected
by that edge. The development of random walks initiated in
different geometric regions are compared and analyzed.

\section{Measures of localization}

\vspace{.2in}\noindent The main problem when studying the
localization of the eigenmodes is to define a suitable notion of
localization. What do we mean when talking about localization of
an eigenfunction? A function on an unbounded domain is said to be
localized if it has compact support or if it is exponentially
localized but a generally accepted notion of localization on
bounded domains is still missing. In fact, we note that due to the
ellipticity of the underlying equation,  there cannot exist any
Dirichlet eigenfunction that vanishes on a nonempty open subset of
the domain \cite{Ev}.

\subsection{The square and triangular Koch prefractals}

\noindent The square and triangular prefractals are defined
recursively by repeatedly applying the corresponding generators.
It was proved that the sequence of square and triangular
prefractals converge to the corresponding fractals and that the
sequences of the individual eigenvalues of the prefractals
converge to the eigenvalues of the actual fractals. Thus, the
study of the eigenfunctions on these approximations should give a
good insight into the behavior of the eigenfunctions on the true
fractals.

\begin{figure}[h]
\label{fig:SQ_SNOW_1_1} \caption{Square Snowflake
Prefractal}\vspace{.3in} \ \ \ \ \ \ \ \ \ \ \ \ \ \ \ \ \ \ \ \ \
\ \ \ \ \ \ \ \ \ \ \ \ \ \ \ \ \ \ \ \ \ \ \ \ \scalebox{.35}
{\includegraphics{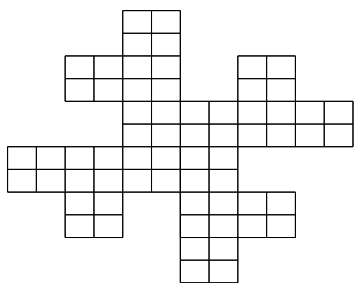}}\label{fig:SQ_SNOW_1_1}
\end{figure}

\begin{figure}[h]
\label{fig:T_SNOW_1_1} \caption{Triadic Snowflake
Prefractal}\vspace{.3in}
 \ \ \ \ \ \ \ \ \ \ \ \ \ \ \ \ \ \ \ \ \ \ \ \ \ \ \ \ \ \ \ \
\ \ \ \ \ \ \ \ \ \ \ \ \ \scalebox{.35}
{\includegraphics{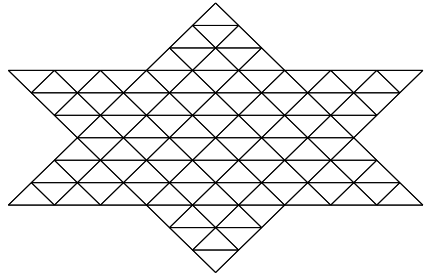}} \label{fig:T_SNOW_1_1}
\end{figure}

\subsection{Localization measures}

\noindent Some eigenmodes on the square Koch prefractal show
confinement to small regions of the membrane. On triadic Koch
prefractals, however, no such `localization' is observed. Sapoval
attributed the experimental localization he observed to the
effects of
\begin{enumerate}
\item Damping \item The existence of `narrow paths' in the
geometry of the structure.
\end{enumerate}

\begin{figure}[h]
\label{fig:LOC_EIG_27} \caption{Eigenfunction 27} \ \ \ \ \ \ \ \
\ \ \ \ \
 \scalebox{1.4} {\includegraphics{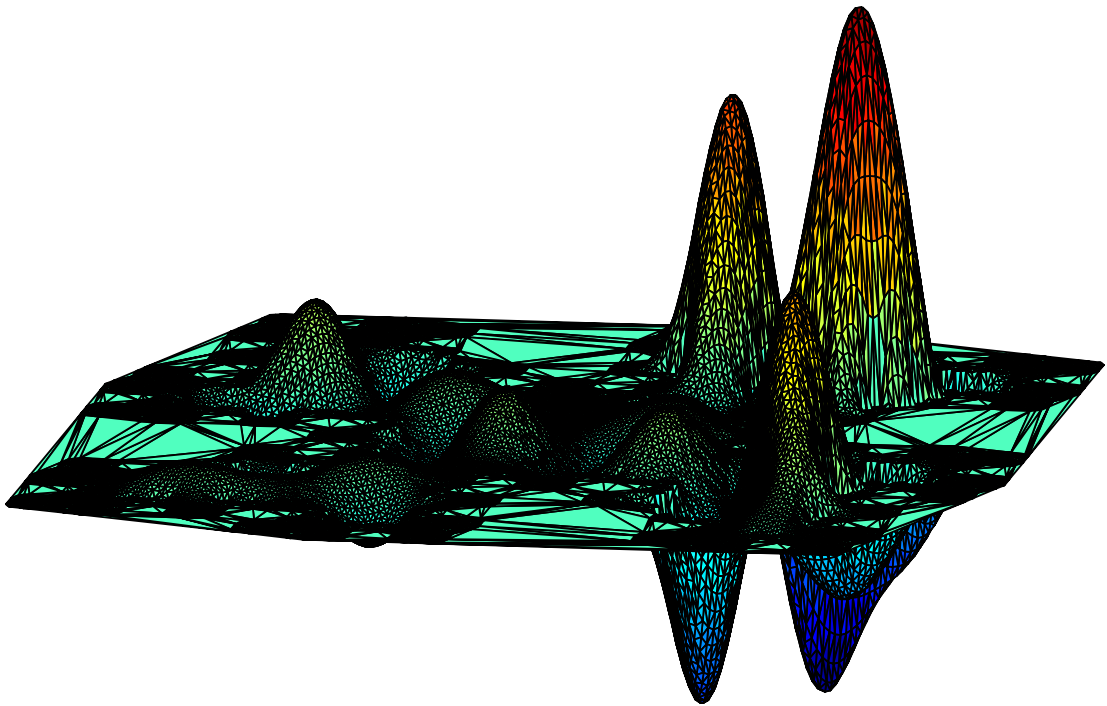}}
 \label{fig:LOC_EIG_27}
\end{figure}

\newpage

\noindent The following quantities were introduced in the
literature in an attempt to qualitatively describe localization of
functions ($\Psi$) on bounded domains:

\vspace{.2in}\textbf{\textit{1) Participation Ratio}} $V_p(\Psi)$
\cite{We}\\
\noindent The Participation Ratio is defined as
$PR(\Psi)=\frac{1}{\int_{D} {|\Psi (x)|}^4dx}.$ It can be regarded
as the effective number of nodes participating in the function
$\Psi$ with significant weight.

\vspace{.2in} \textbf{\textit{2) Participation Volume
$V_p(\Psi)$}} \cite{BeDe}\\
\noindent The Participation Volume is a measure of the effective
number of `atoms' of an `ensemble' participating in a vibration on
a domain $D$: $V_P(\Psi)=\frac{(\int_{D} {({|\Psi
(x)|}^2dx)}^2}{\int_{D} {|\Psi (x)|}^4dx}$. $V_P(\Psi)/V,$ where V
is the volume of the domain $D$, can be used to distinguish
between `\textbf{extended}' (non-localized states) and
`\textbf{exponentially localized states}' (the amplitude decays
exponentially for sufficiently large distances from some central
point).

\vspace{.2in}\noindent In this paper we use the Participation
Ratio in conjunction with suitable notions of diameter and entropy
to investigate localization phenomena on the square and triadic
Snowflake domains.

\section{The Geometry of the Square and Triadic Snowflake Domains}

\noindent How can we describe the existence of narrow channels in
the geometry of the square Snowflake? We will attempt to classify
different regions by defining the following six distances for an
interior grid point $(x,y)$:

\vspace{.2in}\noindent Let $dh_+(x,y)$, $dh_-(x,y)$ be the right
and left horizontal distances of $(x,y)$ from the boundary and
$dv_+$, $dv_-$ be the up and down vertical distances of $(x,y)$
from the boundary. Now define

\begin{itemize}
\item $dh(x,y)=dh_+(x,y)+dh_-(x,y)$

\item $dv(x,y)=dv_+(x,y)+dv_-(x,y)$.
\end{itemize}

\vspace{.2in}\noindent Using $dh(x,y)$ and $dv(x,y)$, one can
distinguish between three geometrically different regions of the
domain in question:
\begin{enumerate}
\item \textbf{`The Grottos'}: Regions confined in both directions. \\
For example, some gridpoints close to the boundary have $dv\approx
dh$ taking small values. \item \textbf{`The Canyons'}: Narrow channels. \\
Points lying in those regions have either $dv\ll dh$ or $dh\ll
dv.$
\item \textbf{`The Prairies'}: Relatively unconfined regions. \\
Some points in the center of the domain, for example, have $dv,
dh\approx 1/2$ and hence should be considered as lying in an
unconfined region.
\end{enumerate}
\vspace{.2in}\noindent\textbf{Remark}: Note that points in region
one as well as points in region three have $dv\approx dh$. Since
the ultimate goal is to define a measure distinguishing between
the three different geometric regions, it is clear that plainly
using $dh(x,y)$ and $dv(x,y)$ will not do.

\vspace{.2in}\noindent We thus define the following three
quantities:

\begin{enumerate}
\item $\partial _{min}(x,y)=\min {\{dh_+(x,y), dh_-(x,y),
dv_+(x,y), dv_-(x,y)\}}$ \item $\partial _{rat}(x,y)=\min
{\{\frac{dh(x,y)}{dv(x,y)}, \frac{dv(x,y)}{dh(x,y)}\} }$  \item
$\partial _{rat}^{min}(x,y)=\partial _{rat}(x,y)*\partial
_{min}(x,y)$.
\end{enumerate}

\vspace{.2in}\noindent Note that for interior grid points in
$SQ(L,R)$, $(L>0)$ with grid size
$H_{SQ(L,R)}={(\frac{1}{2})}^{2L+R}$, we have
\begin{enumerate}
 \item $\partial
_{min}\in [H_{SQ(L,R)},1/4]$\item  $\partial _{rat}\in
[2H_{SQ(L,R)},1]$\item $\partial _{rat}^{min}\in
[2(H_{SQ(L,R)})^2,1/4]$.
\end{enumerate}

\vspace{.2in}\noindent For interior grid points in $T(L,R)$, $(L>0)$
with grid size $H_{T(L,R)}={(\frac{1}{3})}^{L+R}$, we have
\begin{enumerate}
 \item $\partial
_{min}\in [\frac{\sqrt(3)}{2}H_{T(L,R)},\frac{1}{\sqrt{3}}]$\item
$\partial _{rat}\in [\frac{2}{\sqrt{3}}H_{T(L,R)},1]$\item $\partial
_{rat}^{min}\in [(H_{T(L,R)})^2,\frac{1}{\sqrt{3}}]$.
\end{enumerate}

\noindent Let us consider the properties of $\partial_{rat}$ and
$\partial_{rat}^{min}$ in the three different regions:

\noindent\underline{\textbf{Properties of $\partial _{rat}$}}

\begin{itemize}
\item if \textbf{$(x,y)$ lies in a grotto},\\
then $dv \approx dh$ and $\partial _{rat}(x,y)\approx 1.$
\item if \textbf{$(x,y)$ lies in a canyon},\\
then $dv \ll dh$ or $dh\ll dv$ and $\partial _{rat}(x,y)$ is
small.\item if \textbf{$(x,y)$ lies in a prairie},\\
then $dv \approx dh$ and $\partial _{rat}(x,y)\approx 1.$
\end{itemize}


\vspace{.2in}\noindent\underline{\textbf{Properties of $\partial
_{rat}^{min}$}}\\ \noindent Let us assume that the fractal and
refinement levels are not both zero.

\begin{itemize}
\item if \textbf{$(x,y)$ lies in a grotto},\\
then $\partial _{rat}^{min}(x,y)\ \
(\approx{\partial_{min}(x,y)})$\ \ is small.
\item if \textbf{$(x,y)$ lies in a canyon},\\
then $\partial _{rat}^{min}(x,y)\ \ (\approx
{(\partial_{min}(x,y))}^2 )$ is very small.
\item if \textbf{$(x,y)$ lies in a prairie},\\
then $\partial _{rat}^{min}(x,y)\ \ (\approx \partial_{min}(x,y))$
is relatively large.
\end{itemize}

\newpage\noindent Let us now compare the features of $\partial _{rat}^{min}(x,y)$ on
the two domains:

\vspace{.3in}

\begin{tabular}{|c|}
\hline Table 1: Comparison of $\partial _{rat}^{min}(x,y)$ at
$(L,R)=(2,2)$\label{table1}\\ \hline
\end{tabular}

\begin{tabular}{|c|c|c|}
\hline
$\partial _{rat}^{min}(x,y)$ & Square Snowflake & Triadic Snowflake \\
\hline
Maximum & 0.25 & 0.4208 \\
\hline
Minimum & 0.001 & 0.0002 \\
\hline Minimum/Maximum & $4*10^{-3}$ & $4.75*10^{-4}$ \\
\hline Median & 0.0156 & 0.0216\\
\hline Mean & 0.0332 & 0.0724\\
\hline Standard Deviation & 0.0467 & 0.0967 \\
\hline Median/Maximum & 0.0624 & 0.0513 \\
\hline Mean/Maximum & 0.1328 & 0.172 \\
\hline
\end{tabular}

\vspace{.4in}\noindent Note from table 1 that the maximum of
$\partial_{rat}^{min}$ on the square Snowflake is $\frac{1}{4}$
while it is $\frac{1}{\sqrt{3}}$ on the triadic Snowflake. The
minimum value of $\partial_{rat}^{min}$ on the interior of the
triadic domains is smaller than that on the square domain, while
the mean value of $\partial_{rat}^{min}$ is smaller on the square
Snowflake.

\vspace{.1 in}\noindent We find that in both domains, the mean
value of $\partial _{rat}^{min}$ lies closer to the minimum value
of $\partial _{rat}^{min}$ than to the maximum value. Together
with the low medians (only $6.24\%$ of the maximum on the square,
$5.13\%$ of the maximum on the triadic domain), this means that a
lot of points on either domain have a relatively small $\partial
_{rat}^{min}$ value. The relative mean value of $\partial
_{rat}^{min}$ is higher on the triadic Snowflake due to the
triangular geometry. At refinement and fractal level equal to two,
we find that the difference between the maximum and the mean value
lies at about $0.2168$ for the square Snowflake and at about
$0.4206$ for the triadic Snowflake. The differences between
minimum value and mean value lie at about $0.0322$ and $0.0722$,
respectively.

\newpage

\begin{figure}[h]
\label{fig:SQ_dist_distr} \caption{Square Snowflake
$d_{rat}^{min}$ distribution}\vspace{.2in}\ \ \ \ \ \ \ \ \ \ \ \
\ \ \ \ \ \ \ \ \ \ \ \ \ \ \  \scalebox{.55}
 {\includegraphics{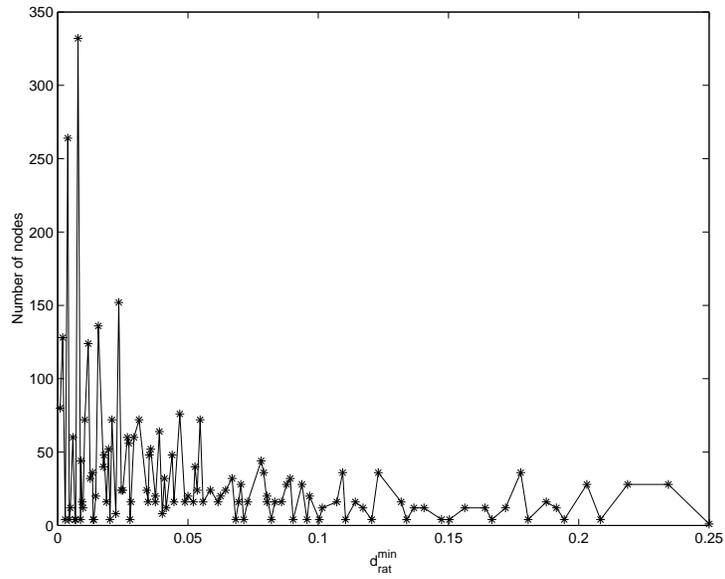}}
\end{figure}

\begin{figure}[h]
\label{fig:D_RAT_MIN_2_2} \caption{Square Snowflake
$d_{rat}^{min}$ surface plot }\vspace{.2in}\ \ \ \ \ \ \ \ \ \ \ \
\ \ \ \ \ \ \ \ \ \ \ \ \ \ \ \ \  \scalebox{1.1}
{\includegraphics{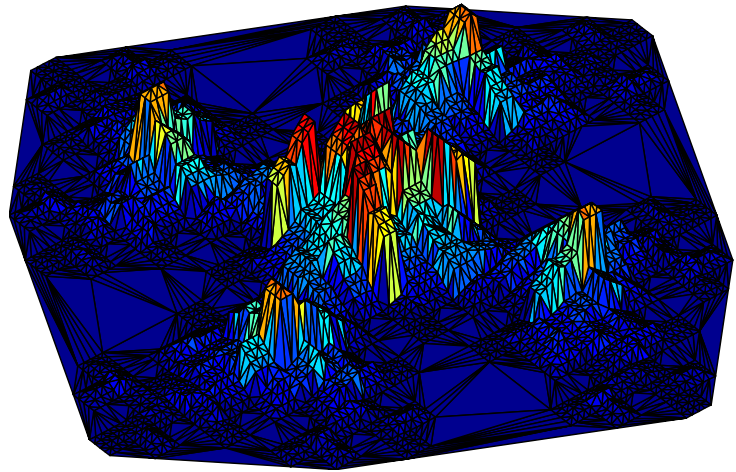}}
\end{figure}

\newpage

\begin{figure}[h]
\label{fig:T_dist_distr} \caption{Triadic Snowflake
$d_{rat}^{min}$ distribution}\vspace{.15in}\ \ \ \ \ \ \ \ \ \ \ \
\ \ \ \ \ \ \ \ \ \ \ \   \scalebox{.55}
{\includegraphics{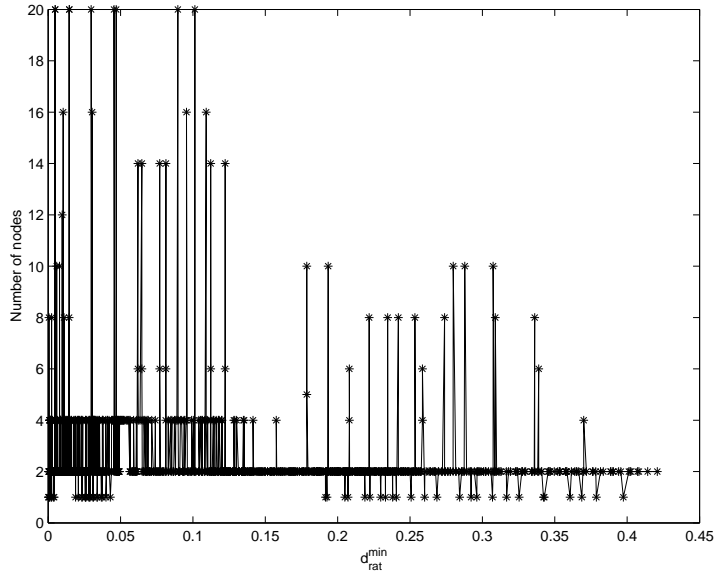}}
\end{figure}

\begin{figure}[h]
\label{fig:fig:T_D_RAT_MIN_2_2} \caption{Triadic Snowflake
$d_{rat}^{min}$ surface plot}\vspace{.2in}\ \ \ \ \ \ \ \ \ \ \ \
\ \ \ \ \ \ \ \ \ \ \ \ \scalebox{.9}
{\includegraphics{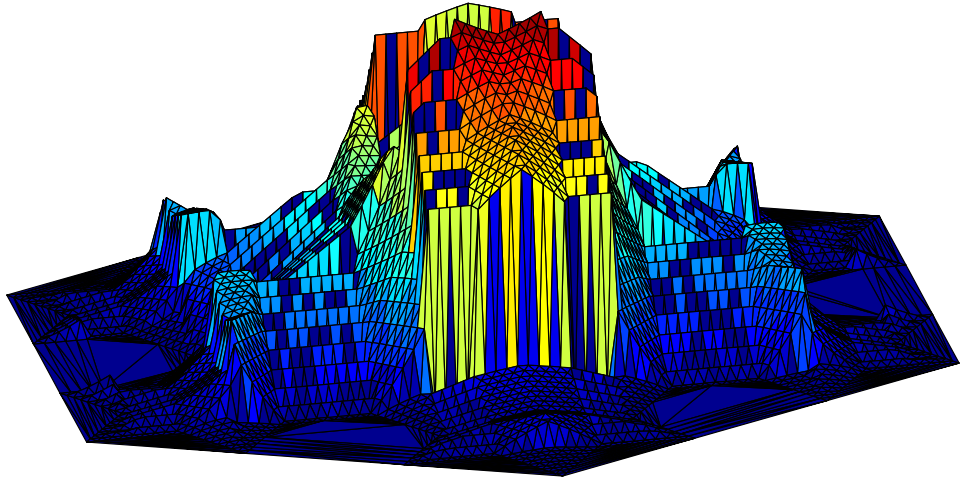}}
\end{figure}

\section{Position Entropy}

\noindent In this section, we will use $\partial_{rat}^{min}$,
defined in the previous section, to define a probability
distribution $P_{x,y}(\partial)$ on each node of our
discretization. This probability distribution is used to define a
position entropy $S(x,y)$ at each vertex.

\vspace{.2in}\noindent For $(x,y)$, a grid point of the
discretization of the snowflake domain, let

\begin{equation}
P_{(x,y)}(\partial):=percentage \quad of\quad neighbours\quad
of\quad node\quad (x,y)\quad with\quad
\partial_{rat}^{min}=\partial .
\end{equation}

\noindent Then $\sum_{\partial}P_{(x,y)}(\partial)=1$ for each
node $(x,y)$ and hence, P defines a probability distribution. Now
let

\begin{equation}
S(x,y):=-\sum_{\partial}P_{(x,y)}(\partial)\ln{(P_{(x,y)}(\partial))}.
\end{equation}

\noindent This entropy gives information about the disorder in the
values of $\partial_{rat}^{min}$ at neighboring points of the
vertex in question. The more variation in  $\partial_{rat}^{min}$,
the larger the value of $S(x,y)$ will be.

\vspace{.2in}

\begin{tabular}{|c|}
\hline Table 2: Comparison of Position Entropy at $(L,R)=(2,2)$ \label{table2}\\
\hline
\end{tabular}

\vspace{.2in}

\begin{tabular}{|c|c|c|}
\hline
Position Entropy & Square Snowflake & Triadic Snowflake \\
\hline Maximum &  0.0157 &   0.0238 \\
\hline Minimum & $2.63*10^{-14}$ & $1.28*10^{-10}$ \\
\hline Minimum/Maximum & $1.675*10^{-12}$ & $5.378*10^{-9}$ \\
\hline Median & 0.0085  & 0.0114 \\
\hline Mean & 0.0089 & 0.0114  \\
\hline Standard Deviation &0.0041& 0.0059\\
\hline Median/Maximum & 0.5414 & 0.47899\\
\hline Mean/Maximum & 0.5668 & 0.47899 \\
\hline
\end{tabular}

\newpage\noindent Comparing our two domains at
$(L,R)=(2,2)$, we find from table 2 that the mean and median
values are the same on the triadic Snowflake and close to each
other on the square Snowflake, and in both domains, they lie
closer to the maximum value of $S$ than to the minimum value. This
means that most of the points on either domain have a relatively
large entropy values, i.e., the values of $\partial_{rat}^{min}$
in the neighborhood of those points are not evenly distributed.
The relative median and mean of $S$ are higher on the square
Snowflake. At refinement and fractal level equal to $2$, we find
that the difference between the maximum and the mean value lies at
about $0.28$ for the square Snowflake and at about $0.38$ for the
triadic Snowflake. The differences between minimum value and mean
value lie at about $0.83$ and $0.95$, respectively.

\vspace{.4in}\noindent In summary, we can say that the relative
mean and median of the entropy found on the square domain are
significantly larger ( about 10$\%$). The $\partial _{rat}^{min}$
characteristics of both domains are similar but we find a slightly
higher (1$\%$) relative median and a lower mean ( about 4$\%$) on
the square Snowflake.

\vspace{.4in}\noindent We refer to figures 9--12 for the
distribution and the surface plot of the square and triadic
Snowflake domains.

\newpage

\begin{figure}[h]
\label{fig:SQ_S_POS_distr} \caption{Square Snowflake Position
Entropy distribution}\vspace{.1in} \ \ \ \ \ \  \ \ \ \ \ \ \ \ \
\ \ \ \ \ \ \ \ \ \ \ \ \ \scalebox{.55}
{\includegraphics{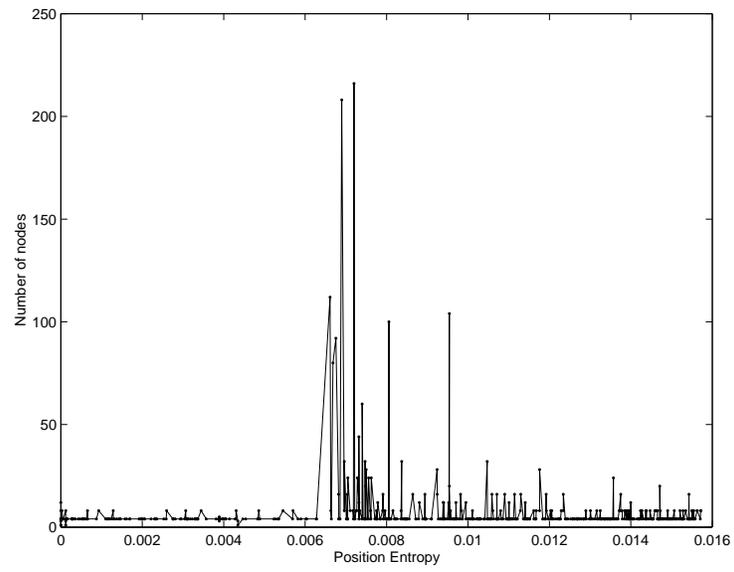}}
\end{figure}

\begin{figure}[h]
\label{fig:_S_POS_2_2} \caption{Square Snowflake Position Entropy
surface plot}\vspace{.2in}\ \ \ \ \ \ \ \ \ \ \ \ \ \ \ \ \ \ \ \
\ \ \ \ \ \ \ \ \ \ \
\scalebox{1}{\includegraphics{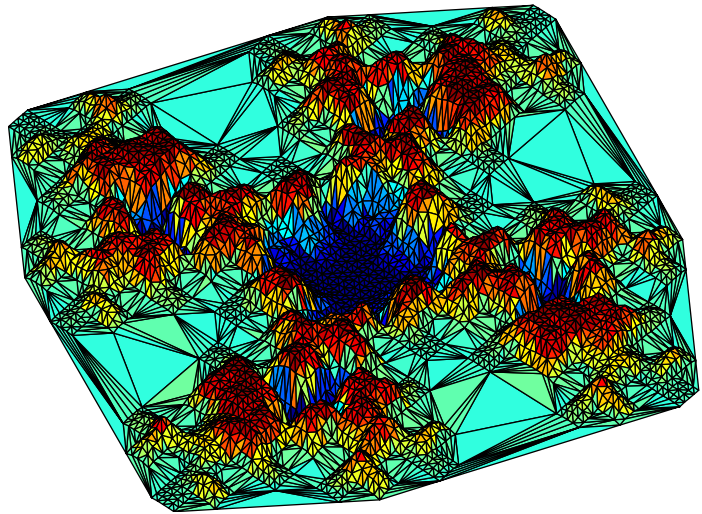}}
\end{figure}

\newpage

\begin{figure}[h]
\label{fig:T_S_POS_distr} \caption{Triadic Snowflake Position
Entropy distribution}\vspace{.1in}  \ \ \ \ \ \  \ \ \ \ \ \ \ \ \
\ \ \ \ \ \ \ \ \ \ \ \ \ \scalebox{.55}
{\includegraphics{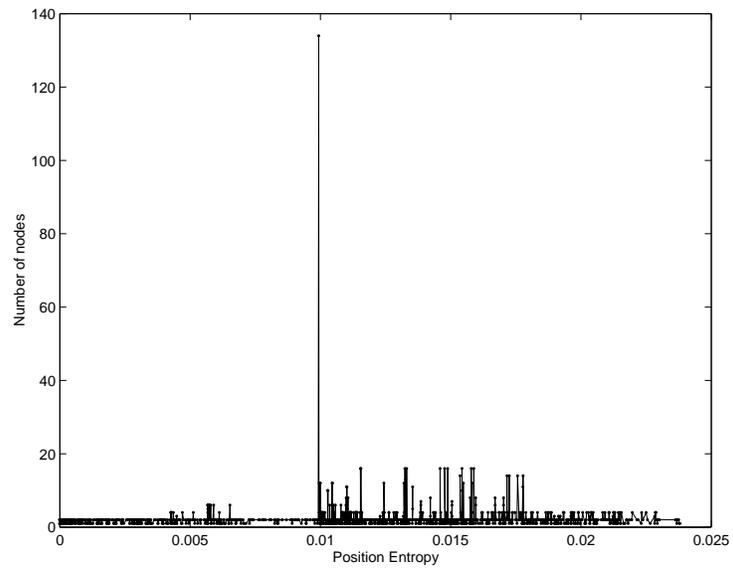}}
\end{figure}

\vspace{.2in}

\begin{figure}[h]
\label{fig:T_S_POS_2_2} \caption{Triadic Snowflake Position
Entropy surface plot}\vspace{.2in}\ \ \ \ \ \ \ \ \ \ \ \ \ \ \ \
\ \ \ \ \ \ \scalebox{1}{\includegraphics{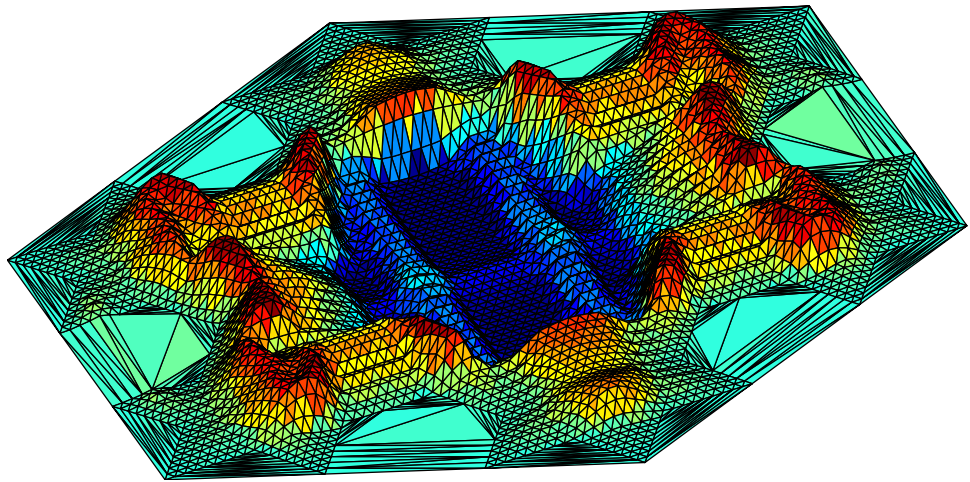}}
\end{figure}

\section{Random Walks on Discretizations of Snowflake Domains}

\noindent In this section, we will study diffusion (random walks)
on discretization graphs of the Snowflake domains. This work is
inspired by I. Simonsen et. al. \cite{SiErMaSn}, \cite{Si}. In
their papers, these authors use diffusion on networks to derive a
master equation whose analysis reveals information about the
large-scale structure of the network. By performing random walks
governed by the underlying master equation on the two Snowflake
domains, we hope to gain useful insights into their geometric
properties and an understanding of why localization is favored by
the Square domain. This will provide a new way to investigate
localization for domains with fractal boundaries. This study
relies in part on the notions of distance and entropy introduced
earlier in this paper (see Sections 3 and 4).

\vspace{.1in}\noindent The nodes of the discretizations represent
the vertices and the edges are given by the neighbor relations.
The path taken by our random walkers will be influenced by the
geometric features of the domain:

\vspace{.1 in} \noindent The edge $e_{ij}$, going from vertex
$v_i$ to vertex $v_j$ is assigned a weight inversely proportional
to $1-\partial_{rat}^{min}(j)$, where $\partial_{rat}^{min}$ is
normalized such that it takes values in $[0,1]$ . Hence, edges
where node $v_j$ lies in a Grotto have large weights assigned,
edges where node $v_j$ lies in a Canyon have even larger weights
(close to 1) assigned, while edges with node $v_j$ lying in a
Prairie have relatively small weights assigned, i.e., the weights
are assigned to the edges in such a way that a random walker is
strongly encouraged to move to regions with distances that
describe narrow channels, then to regions confined in both
directions, and lastly to unconfined regions. If localization is
caused by the existence of narrow channels, as predicted by
Sapoval, we expect to see some kind of congregation of the walkers
in those areas of the domain that are characterized by
predominantly large weights.

\vspace{.2 in} \noindent We start our random walk by placing a
large number of random walkers onto the vertices of the network.
At each time step, these walkers are allowed to move between
adjacent vertices. The edge, out of the possible outgoing ones, a
walker chooses to move along is picked at random with a
probability equal to the weight assigned to this edge.

\vspace{.2in}\noindent We define the adjacency matrix $A$ as
follows:

\begin{equation}
A_{ij}= \Bigg\{
\begin{array}{cc}
1, & \mbox{if}\quad nodes\quad v_i\quad and\quad v_i\quad are\quad neighbors\nonumber\\
0, & else\nonumber .\\
\end{array}
\end{equation}

\noindent Note that the weight matrix
$W=\{w_{ij}\}_{i,j}=\{A_{ij}(1-\partial_{rat}^{min}(j))\}$ is
non-symmetric.

\subsection{The Master Equation}

\noindent Let $nv$ be the number of vertices in the discretization
and $N$ be the number of participating random walkers. Let the
number of walkers situated at node $v_i$ at time $t$ be denoted by
$N_i(t)$. Then the fraction of walkers at this node at time $t$ is

\begin{equation}
\eta_i(t)=\frac{N_i(t)}{N}.
\end{equation}

\noindent Since the total number of random walkers is preserved at
each time, we have $\sum_{i}\eta_i(t)=1$.

\vspace{.1in}\noindent The change in the walker density of a
vertex $v_i$ during one time step equals the difference between
the relative number of walkers entering and leaving the same
vertex over the time interval. In mathematical terms we can write

\begin{equation}
\eta_i(t+1)=\eta_i(t)+J_i^-(t)-J_i^+(t), \label{41a}
\end{equation}

\noindent where $J_i^{\pm}(t)$ denote the relative number of
walkers entering (-) and leaving (+) vertex $v_i$.

\vspace{.1in}\noindent The value of $J_i^+(t)$ depends on the
\textbf{strength} $St_i$ of node $v_i$ i.e.,  the total number of
outgoing weights $St_i=\sum_{j\in Neigh_i}w_{ij}$, where $Neigh_i$
denotes the set of neighbors of node $v_i$. Similarly, the value
of $J_i^-(t)$ depends on the total number of incoming weights
$\sum_{j\in Neigh_i}w_{ji}$. The fraction of outgoing walkers from
vertex $v_i$ (a current) per unit weight is thus

\begin{equation}
c_i(t)=\frac{\eta_i(t)} {\sum_{j\in Neigh_i}w_{ij}}.
\end{equation}

\noindent Hence, $c_i(t)$ is the weighted walker density per link
of node $v_i$ at time $t$.

\vspace{.1in}\noindent The edge current on the directed edge from
vertex $v_i$ to vertex $v_j$ is then given by

\begin{equation}
C_{ij}(t)=w_{ij}c_i(t)=w_{ij}\frac{\eta_i(t)}{\sum_{k\in
Neigh_i}w_{ik}}.
\end{equation}

\noindent Notice that the term $$\frac{w_{ij}}{\sum_{k\in
Neigh_i}w_{ik}}$$ is the probability of a walker moving from
vertex $v_i$ to vertex $v_j$.

\vspace{.1in}\noindent The relative number $J^+_i$ of outgoing
walkers from node $v_i$  at time $t$ is given by

\begin{equation}
J^+_i(t)=\sum_{j\in Neigh_i}C_{ij}
\end{equation}

\noindent and the relative number $J^-_i$ of incoming walkers to
node $v_i$ at time $t$ is given by

\begin{equation}
J^-_i(t)=\sum_{j\in Neigh_i}C_{ji}.
\end{equation}

\noindent It is easy to verify that $J^+_i(t)=\eta_i(t)$. This
expresses the fact that all walkers leave their respective
vertices at each time step, i.e., no walker stays at the same
vertex for more than one time step.

\vspace{.1in}\noindent  Denoting
$\partial_t\eta_i(t)=\eta_i(t-1)-\eta_i(t)$ and substituting the
expressions for $ J^+_i(t)$ and $J^-_i(t)$ into equation
(\ref{41a}), we find

\begin{eqnarray}
\partial_t\eta_i(t)&=&\sum_{j\in Neigh_i}C_{ji}-\eta_i(t)\nonumber\\
&=& \sum_{j\in Neigh_i}\frac{w_{ji}\eta_j(t)}{\sum_{k\in
Neigh_j}w_{jk}}-\eta_i(t)\nonumber\\
&=&\sum_{j\in Neigh_i}T_{ij}\eta_j(t)-\eta_i(t),
\end{eqnarray}

\noindent where $$T_{ij}:=\frac{w_{ji}}{\sum_{k\in
Neigh_j}w_{jk}}.$$  \noindent This equation can be written in
matrix form, as follows:

\begin{equation}
\mathbf{\eta}(t+1)=\mathbf{T\eta}(t)\label{41b}
\end{equation}

\noindent Equation (\ref{41b}) is known as the master equation for
the random walk process of the underlying network. Further,
$\mathbf{T}=\{T_{ij}\}$ is called the transfer matrix. It
`transfers' the walker distribution one step ahead and therefore
can be thought of as a time propagator for the process. For some
arbitrarily chosen initial state $\mathbf{\eta}(0)$, the time
development can be obtained by iterations of equation (\ref{41b}),
with the result
$\mathbf{\eta}(t)=\mathbf{T}^{(t)}\mathbf{\eta}(0)$. Hence, the
eigenvalue spectrum of \textbf{T} controls the time evolution of
the diffusive process.

\vspace{.1in}\noindent Due to the fact that our weight matrix is
not symmetric, the matrix \textbf{T} also fails to be symmetric.
However, $\mathbf{T}$ is similar to the symmetric matrix
$\mathbf{S=KTK^{-1}}$, where

\begin{equation}
K_{ij}:= \Bigg\{
\begin{array}{cc}
\frac{\delta_{ij}}{\sqrt{w_{ij}}*\sqrt{\sum_{k\in Neigh_j}{w_{jk}}}}, & \mbox{if}\quad w_{ij}\neq 0\nonumber\\
1, & else,\nonumber\\
\end{array}
\end{equation}

\noindent and hence
\begin{eqnarray}
{(KTK^{-1})}_{ij}&=& \frac{\sqrt{w_{ii}w_{jj}}}{\sqrt{\sum_{k\in
N_i}w_{ik}\sum_{k\in N_j}w_{jk}}}\nonumber\\
&=&
\frac{\sqrt{(1-\partial_{rat}^{min}(i))(1-\partial_{rat}^{min}(j))}}{\sqrt{[\sum_{k\in
N_i}(1-\partial_{rat}^{min}(k))][\sum_{k\in N_j}(1-\partial_{rat}^{min}(j))]}}.\nonumber\\
\end{eqnarray}

\noindent Therefore, \textbf{T} is guaranteed to have real
eigenvalues and corresponding eigenvectors. It is convenient to
sort the eigenvalues in descending order. Furthermore, as a
consequence of the Perron--Frobenius theory, the eigenvalue
$\lambda_1=1$ is simple and the elements of the corresponding
eigenvector will all have the same sign.

\vspace{.1in}\noindent Physically, the state $\lambda_1=1$
corresponds to the stationary state such that
$\mathbf{\eta}(\infty)\propto\mathbf{\eta^{t}}$, where the
diffusive current flowing from node $v_i$ to node $v_j$ is exactly
balanced by that flowing from $v_j$ to $v_i$. Since
$|\lambda_k|<1$ for $k\neq 1$, all modes corresponding to those
eigenvalues are decaying. $\lambda_k>0$ correspond to
non-oscillating modes while $\lambda_k<0$ correspond to states
where oscillation will take place over time. The large scale
topology of the given network reflects itself in the statistical
properties of the eigenvectors $\eta_k$ \cite{SiErMaSn},
\cite{Si}.

\section{Numerical Work}

\noindent As shown in the previous section, a random walk on the
discretization  of the square and triadic Snowflake domains is
governed by the master equation
$\mathbf{\eta}(t)=\mathbf{T}^{(t)}\mathbf{\eta}(0)$, where
$\mathbf{\eta}(t)$ is the distribution of random walkers at time
$t$. In this section, we will discuss numerical simulations of
these random walks. Initially, we will place all random walkers at
one vertex. By repeatedly applying $T$ to our distribution, we are
able to obtain the new walker distributions for each time step. We
expect that the dynamics of this system is strongly influenced by
the choice of the initial node. Recall that a node lying in a
narrow channel is characterized by small $\partial_{rat}^{min}$
and entropy values. Hence, in accordance with the definition of
the weight matrix $W$ ($w_{ij}=1-\partial_{rat}{min}(j)$), if
narrow channels in the geometry of the square Snowflake exist, we
would hope to see a concentration of walkers in the regions where
$\partial_{rat}^{min}$ and entropy values are low.

\newpage
\vspace{1in}
\begin{figure}[h]
 \caption{Triadic Snowflake: Position
Entropy plus Ratio Distance contour plot} \vspace{.3in} \ \ \ \ \
\ \ \ \ \ \ \ \ \ \ \ \ \ \ \ \ \ \ \  \scalebox{.8}
{\includegraphics{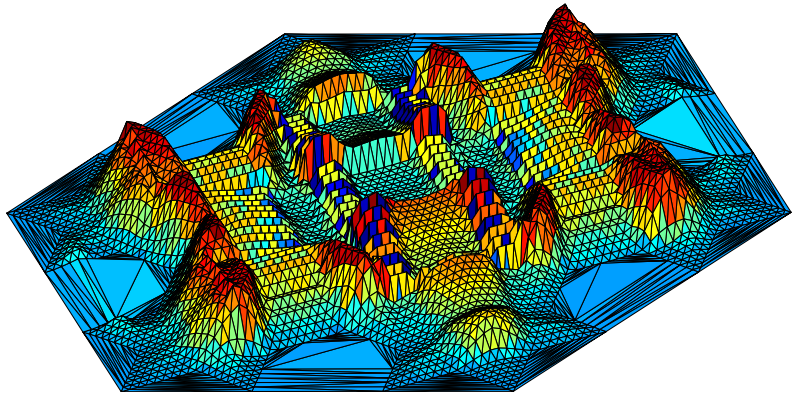}}
\label{fig:T_S_POS_D_RAT_MIN}
\end{figure}

\vspace{1in}\noindent In figure 13, we see a plot of the relative
Minimum Ratio Distances added to the relative Position Entropies
at each grid point of the Triadic Snowflake. As expected, we see
low values close to the boundary. Indeed, random walkers situated
here are expected to stay localized in the region. The center part
of the Snowflake comprises a medium value square on the left and
right of which we find two rectangular discrete low bands and 4
other medium value regions situated north, west, south and east of
the center square. In summary, the triadic Snowflake has a
relatively big open area of low and medium values in its center.
This is not the case for the square Snowflake.

\newpage

\begin{figure}[h]
\label{S_POS_D_RAT_MIN} \caption{Square Snowflake: Position
Entropy plus Ratio Distance contour plot} \vspace{.3in} \ \ \ \ \
\ \ \ \ \ \ \ \ \ \ \ \ \ \ \ \ \ \ \ \ \ \ \ \ \ \ \ \
\scalebox{.8} {\includegraphics{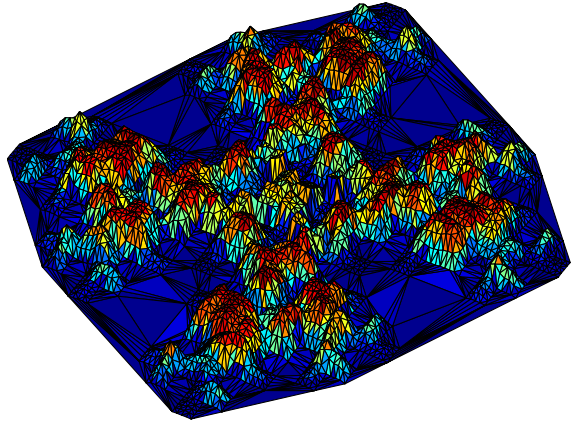}}
\end{figure}

\vspace{1in}\noindent Figure 14
shows a plot of the relative Minimum Ratio Distances added to the
relative Position Entropies at each grid point of the square
Snowflake. Note that in the center we find relatively high values
but, surrounding that center, we find a band of low entropy and
distance values. We see four arms of high entropy/distance
stretching through parts of the domain almost reaching the
boundary. In each of those four arms we can find a small low value
basin. The remaining regions between the arms and near the
boundary  are also characterized by low values. We would expect
that a random walk started in any of the basins would stay
localized, i.e., that the majority of high/low amplitudes can be
found there. If the random walk is initiated in a high altitude
region, we would hope to see no such localization. If this were
observed, it would indeed indicate that localization really is
caused by the existence of narrow channels.

\newpage\noindent To test this theory, we initially place all random
walkers at a point of interest (low entropy/low distance, high
entropy/high distance, low entropy/high distance or high
entropy/low distance) and then watch the random walk evolve by
computing the following quantities at time t:

\begin{itemize}
\item the diameter $DIAM(t)$, the furthest distance from the
initial node traveled by a random walker  during $[0,t]$. \item
the Participation Ratio $PR(t)=\frac{1}{ \sum_i{[\eta_i(t)]}^4}$,
where $\eta_i(t)$ is the percentage of random walkers at node
$v_i$ at time $t$.
\end{itemize}

\noindent We define high amplitude points at time $t$ as points
where $|\eta_i(t)|>(3/4)max(|\eta_i(t)|)$ and compute

\begin{itemize}
\item $D_{ha}^{rel}(t)$, the Relative Minimum Ratio Distance of
high amplitude points at time t, i.e.,
$D_{ha}^{rel}(t)=\frac{mean({(\partial_{rat}^{min})}_{ha}(t))}{max(\partial_{rat}^{min})}$,
where $mean({(\partial_{rat}^{min})}_{ha}(t))$ is the mean of the
Minimum Ratio Distances  of high amplitude points. \item
$S_{ha}^{rel}(t)$, the Relative Position Entropy of high amplitude
points at time t, i.e.,
$S_{ha}^{rel}(t)=\frac{mean(S_{ha}(t))}{max(S)}$, where
$mean(S_{ha}(t))$ is the mean of the entropies of high amplitude
points. \item $DIAM_{ha}^{rel}(t)=\frac{DIAM_{ha}(t)}{DIAM(t)}$,
the Relative Diameter of high amplitude points at time $t$
($DIAM_{ha}(t)$ denotes the diameter of high amplitude points).
\end{itemize}

\vspace{.2in}\noindent We also keep track of where the walk was
started. Since we want to be able to compare values on the two
different domains, we will need to compute relative rather than
actual values: Instead of listing $\partial_{rat}^{min}(1)$ and
$S(1)$, we list
$\partial_{rel}(1)=\frac{\partial_{rat}^{min}(1)}{\max(\partial_{rat}^{min})}$
and $S_{rel}(1)=\frac{S(1)}{max(S)}$.

\vspace{.1in}\noindent The following six figures, numbered 15--20,
show some typical plots obtained from our numerical calculation.
We would like to point out that there is a lot of variation in the
character of the plots: On both domains, the plots of the four
quantities look different for different starting nodes but in
almost all cases, the four quantities approach an asymptotic limit
after about 200 time steps.

\newpage

\begin{figure}[h]
\label{fig:SQ_SIM} \caption{Square Snowflake: Random Walk
Simulation}\vspace{.1in}\ \ \ \ \ \ \ \ \ \ \ \ \ \ \ \ \ \ \ \ \
\ \ \ \ \ \ \ \ \ \ \ \ \ \ \ \ \ \ \ \ \ \ \ \scalebox{.35}
{\includegraphics{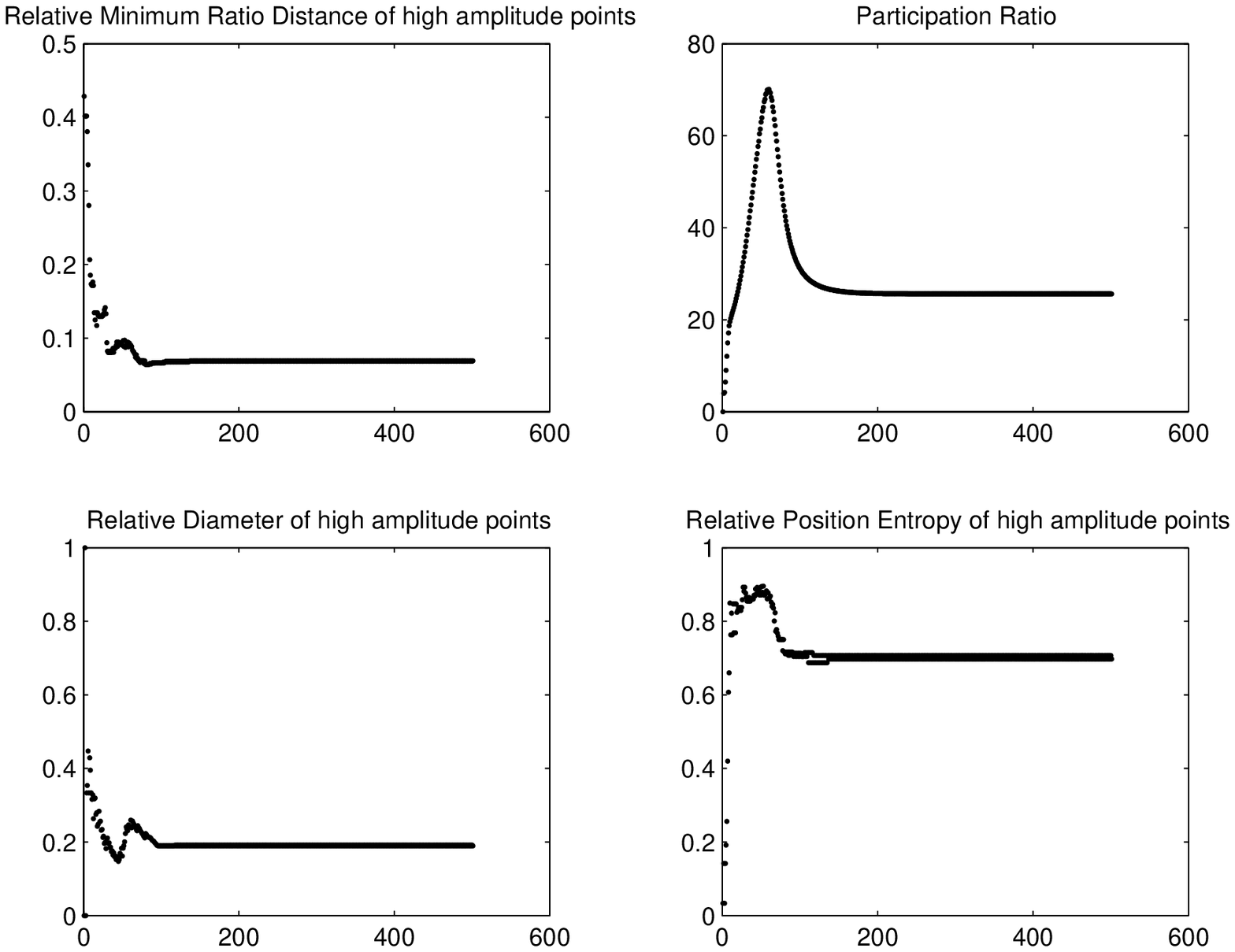}}
\end{figure}

\begin{figure}[h]
\label{fig:SQ_SIM_2} \caption{Square Snowflake: Random Walk
Simulation}\vspace{.1in}\ \ \ \ \ \ \ \ \ \ \ \ \ \ \ \ \ \ \ \ \
\ \ \ \ \ \ \ \ \ \ \ \ \ \ \ \ \ \ \ \ \ \ \ \scalebox{.35}
{\includegraphics{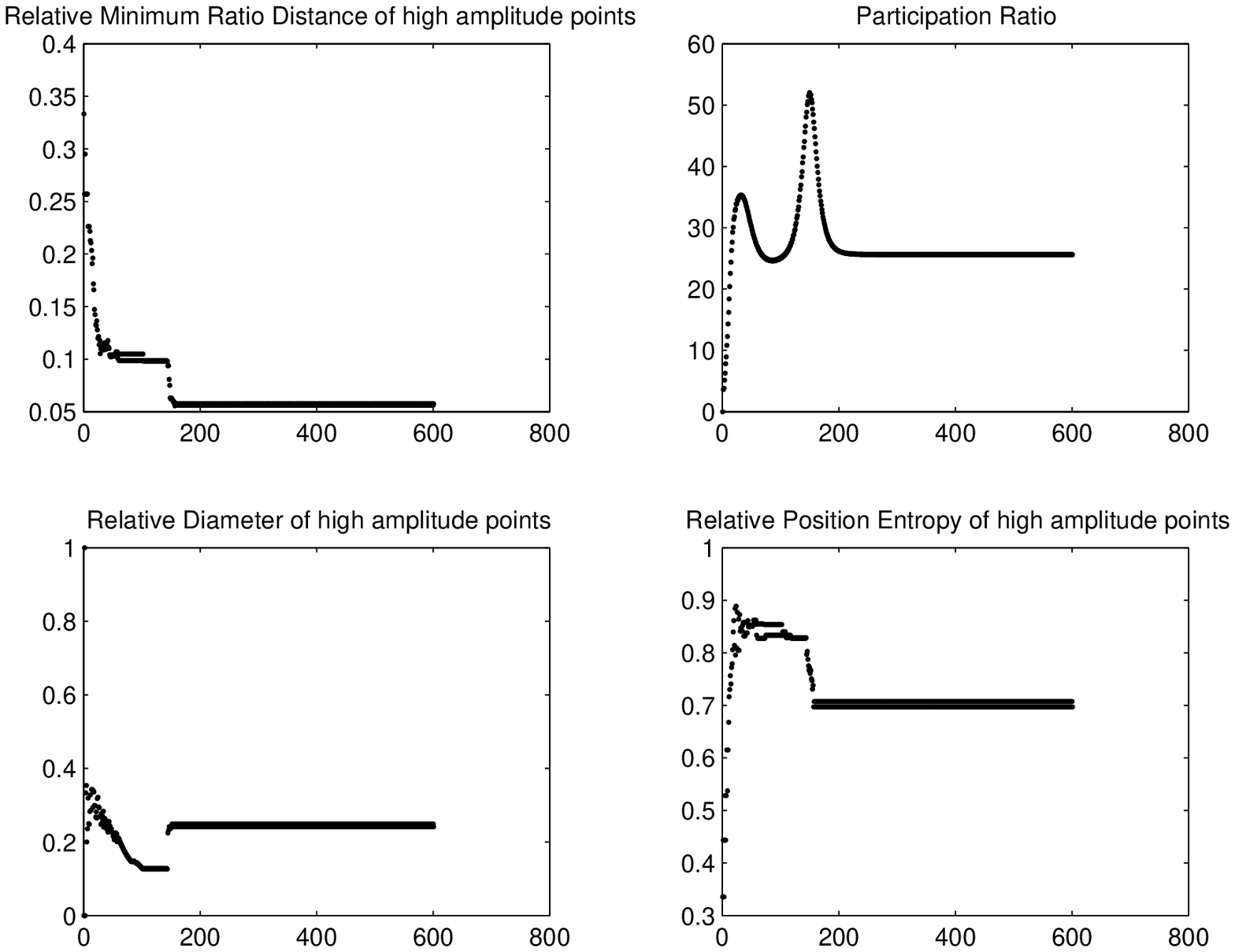}}
\end{figure}

\begin{figure}[h]
\label{fig:DFF_SIM_1} \caption{Square Snowflake: Random Walk
Frame}\ \ \ \ \ \ \ \ \ \ \ \ \ \ \ \ \ \ \ \ \ \ \ \ \ \ \ \ \ \
\ \ \ \ \ \ \ \ \ \ \ \ \ \ \scalebox{.35}
{\includegraphics{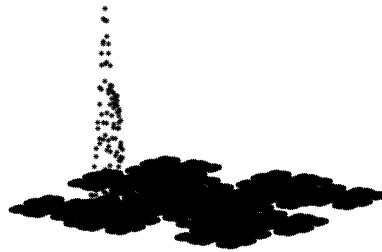}}
\end{figure}

\newpage

\begin{figure}[h]
\label{fig:T_SIM} \caption{Triadic Snowflake: Random Walk
Simulation}\vspace{.1in}\ \ \ \ \ \ \ \ \ \ \ \ \ \ \ \ \ \ \ \ \
\ \ \ \ \ \ \ \ \ \ \ \ \ \ \ \ \ \ \ \ \ \ \ \scalebox{.35}
{\includegraphics{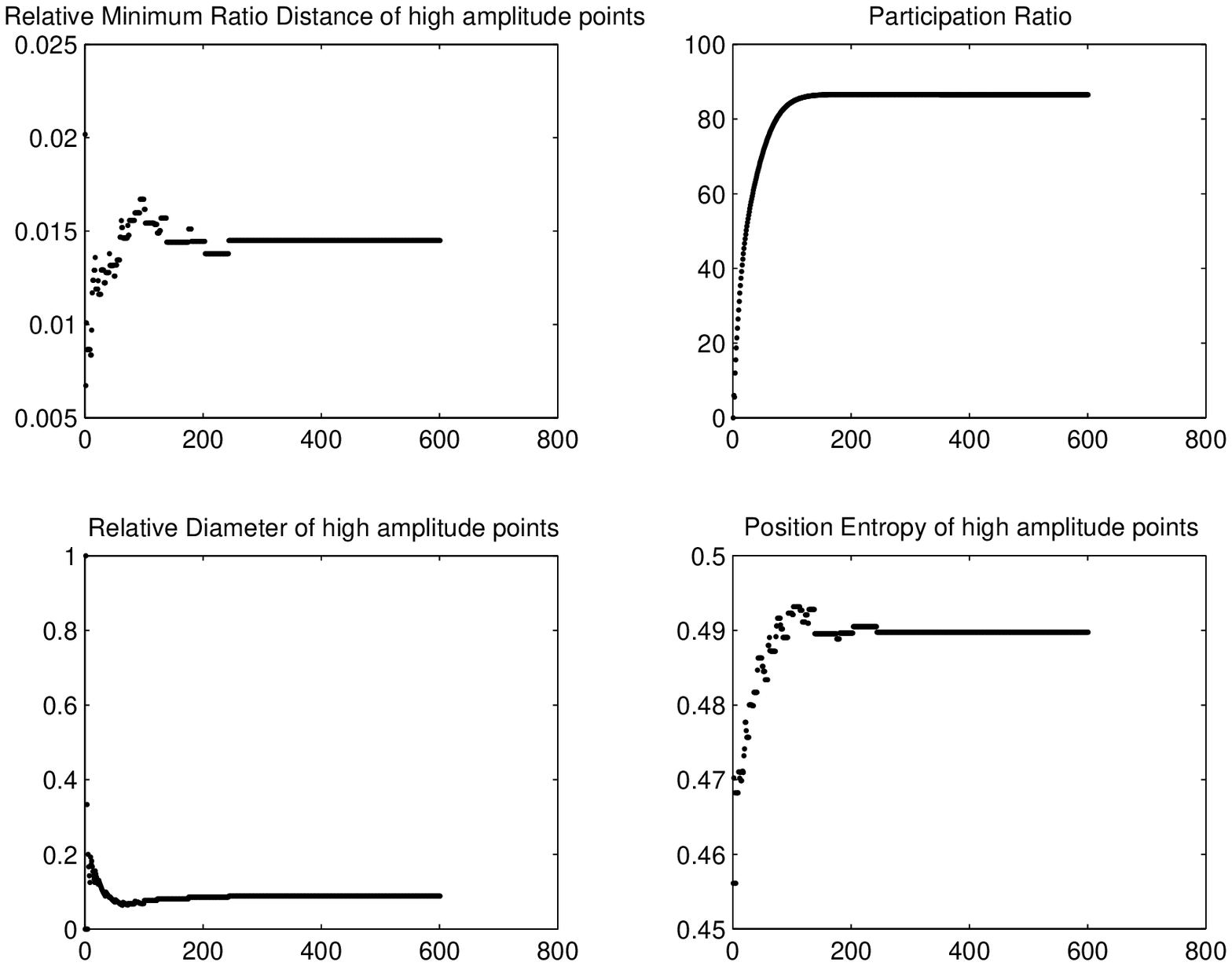}}
\end{figure}

\begin{figure}[h]
 \caption{Triadic Snowflake: Random Walk
Simulation}\vspace{.1in}\ \ \ \ \ \ \ \ \ \ \ \ \ \ \ \ \ \ \ \ \
\ \ \ \ \ \ \ \ \ \ \ \ \ \ \ \ \ \ \ \ \ \ \ \scalebox{.35}
{\includegraphics{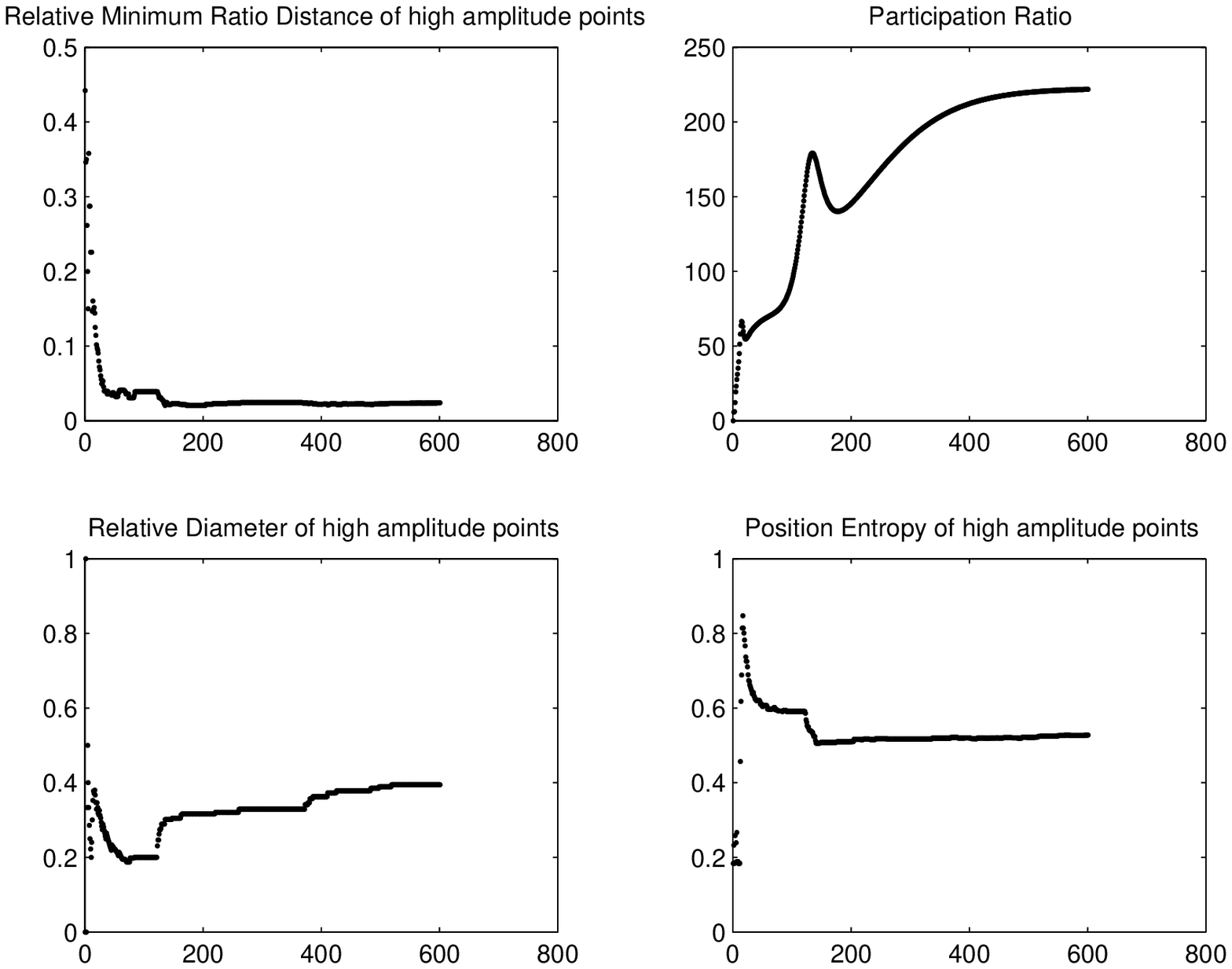}} \label{figT_SIM_2}
\end{figure}

\begin{figure}[h]
\label{fig:T_DFF_SIM_1} \caption{Triadic Snowflake: Random Walk
Frame}\ \ \ \ \ \ \ \ \ \ \ \ \ \ \ \ \ \ \ \ \ \ \ \ \ \ \ \ \ \
\ \ \ \ \ \ \ \ \ \ \ \ \ \ \scalebox{.35}
{\includegraphics{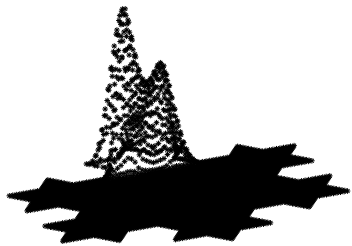}}
\end{figure}

\noindent Figure \ref{figT_SIM_2} shows a simulations that takes
more than 200 time steps for the values to approach the steady
state.

\newpage\noindent Let us look at some random walks started at low
entropy/low $\partial_{rat}^{min}$ nodes on both domains.

\begin{tabular}{|c|}
\hline
Table 3: \textbf{Square Snowflake} ($\partial_{rat}^{min}(1)$ decreases)\label{table3}\\
\hline
\end {tabular}

\vspace{.1in}

\begin{tabular}{|c|c|c|c|c|c|c|}
\hline $\partial_{rel}(1)$ & $S_{rel}(1)$& $PR(\infty)$ &
$DIAM(\infty)$ & $D_{ha}^{rel}(\infty)$ & $S_{ha}^{rel}(\infty)$ &
$DIAM_{ha}^{rel}(\infty)$ \\
\hline 0.4572 & $8.21*10^{-3}$ & 25.6456 & 0.9254 & 0.0559 &
0.707 & 0.5829\\
\hline 0.3216 & 0.12& 25.5995 & 0.9512 & 0.0667 & 0.69 & 0.54\\
\hline 0.3124 & 0.146 & 25.5997 & 0.9071 & 0.0557 & 0.6972 &
0.6225\\
\hline 0.156 & 0.4777 & 25.5996 & 1.3436 & 0.0541 & 0.6972 &
0.3871 \\
\hline $8*10^{-4}$ & 0.4777 & 25.6527 & 1.3295 & 0.0562 & 0.707 &
0.3716\\
\hline $9.7*10^{-4}$ & 0.4777 & 25.6455 & 1.5202 & 0.0493 & 0.7071
& 0.0783\\ \hline
\end{tabular} \\

\noindent The diameter increases with decreasing
$\partial_{rat}^{min}(1)$ values but the relative diameter of high
amplitude points decreases. Note that the entropy of high
amplitude points is lying at around $70\%$ of the maximum value
which is much larger than we would expect, but the minimum ratio
distances are very low ($\approx 5\%$ of the maximum value).

\vspace{.1in}\noindent Also note the large diameter but very low
relative diameter of the high amplitude points in the last row.
The corresponding plots for this simulation are given in figure
21.

\begin{figure}[h]
\label{SQ_L_L_2415} \caption{Square Snowflake: Random Walk
Simulation node 2415}\ \ \ \ \ \ \ \ \ \ \ \ \ \ \ \ \ \ \ \ \ \ \
\ \scalebox{.5} {\includegraphics{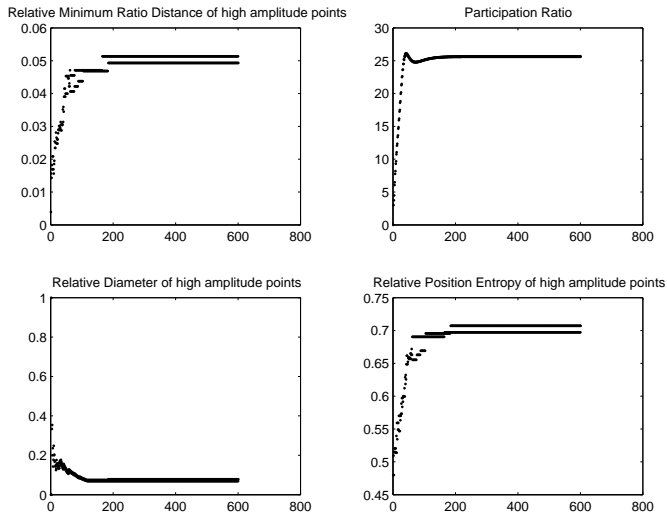}}
\end{figure}

\noindent At a time step of about 200, $D_{ha}^{rel}$ and
$S_{ha}^{rel}$ start jumping between two values; an indication of
localization?

\vspace{.2in}

\begin{tabular}{|c|}
\hline Table 4: \textbf{Square Snowflake}
($S(1)$ decreases)\label{table4}\\
\hline
\end {tabular}

\vspace{.1in}

\begin{tabular}{|c|c|c|c|c|c|c|}
\hline $\partial_{rel}(1)$ & $S_{rel}(1)$ & $PR(\infty)$ &
$DIAM(\infty)$ & $D_{ha}^{rel}(\infty)$ & $S_{ha}^{rel}(\infty)$ &
$DIAM_{ha}^{rel}(\infty)$ \\
\hline 0.314 & 0.1464 & 25.5997 & 0.9071 & 0.0557 & 0.6972 &
0.6225\\
\hline 0.3436 & 0.07 & 25.6456 & 0.888 & 0.0566 & 0.7071 &
0.6483\\
 \hline 0.3572 & $0.058$ & 25.5997 & 1.1761 &
0.0575 &
0.6972 & 0.2941\\
 \hline 0.4572 & $8.21*10^{-3}$ & 25.6456 &
0.9254 & 0.0559 &
0.707 & 0.5829\\
\hline 0.492 & $1.78*10^{-3}$ &  25.6456  & 0.9135 & 0.0559 &
0.7071 & 0.6047\\
\hline 0.492 & $1.52*10^{-3}$ & 25.6005 & 0.9209 & 0.0568 & 0.6972
& 0.5998\\
 \hline
\end{tabular} \\

\noindent The diameter of high amplitude points is large and does
not decrease with decreasing entropy values. Here, row 3 sticks
out with a relatively low $DIAM_{ha}^{rel}(\infty)$ value. The
corresponding plots are given in figure 22.

\begin{figure}[h]
\label{SQ_L_L_824} \caption{Square Snowflake: Random Walk
Simulation node 824}\ \ \ \ \ \ \ \ \ \ \ \ \ \ \ \ \ \ \ \ \ \ \
\ \scalebox{.5} {\includegraphics{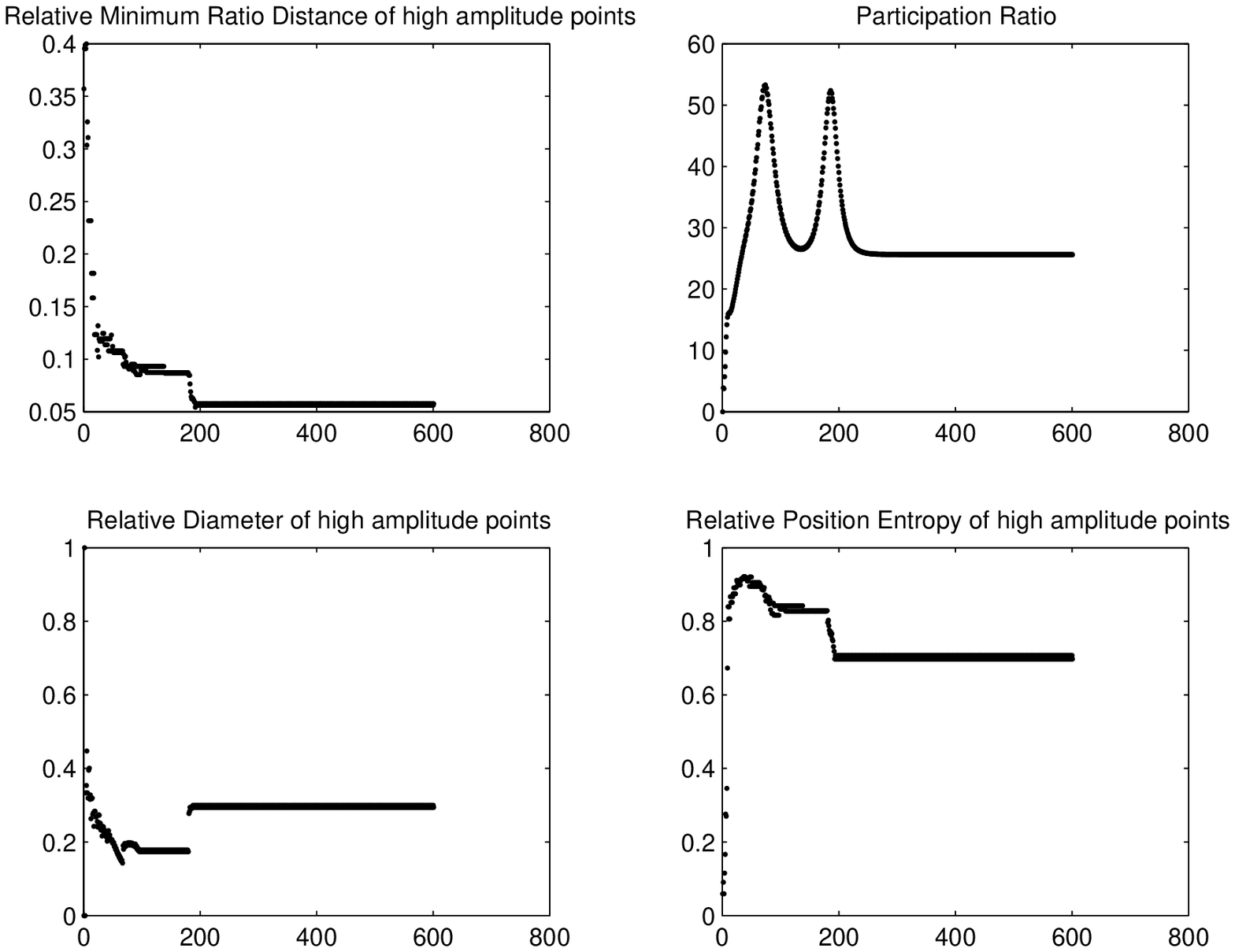}}
\end{figure}

\noindent In both tables, $D_{ha}^{rel}(\infty)$ and
$S_{ha}^{rel}(\infty)$ show little fluctuations.

\newpage

\begin{tabular}{|c|}
\hline
Table 5: \textbf{Triadic Snowflake} ($\partial_{rat}^{min}(1)$ decreases) \label{table5}\\
 \hline

\end {tabular}

\vspace{.1in}

\begin{tabular}{|c|c|c|c|c|c|c|}
\hline $\partial_{rel}(1)$ & $S_{rel}(1)$ & $PR(\infty)$ &
$DIAM(\infty)$ & $D_{ha}^{rel}(\infty)$ & $S_{ha}^{rel}(\infty)$ &
$DIAM_{ha}^{rel}(\infty)$ \\
\hline 0.4776 & 0.063 & 156.8692 & 0.5492 & 0.0153 & 0.4897 &
0.7333\\
\hline 0.3973 & 0.231 & 169.2816 & 0.7011 & 0.0144 &
0.4883 & 0.5033\\
\hline 0.0206 & 0.4579 & 221.154 & 0.9773 & 0.0231 & 0.5247 &
0.1822\\
\hline 0.0163 & 0.4663 & 86.536 & 0.011 & 0.0147 & 0.4897 &
0.1895\\
\hline 0.0617 & 0.4285 & 221.3947 & 0.9558 & 0.023 & 0.5246 &
0.1944\\
\hline $1.13*10^{-3}$ & 0.4201 & 222.6706 & 1.0908 & 0.0308 &
0.5301 & 0.187\\\hline
\end{tabular} \\

\noindent In contrast with what was observed for the square
domain, we see here in table 5 different participation ratios for
different starting nodes. Note that on the triadic domain the
entropy of high amplitude points are lying at only about $50\%$,
as opposed to about $70\%$ on the other domain. As observed on the
square domain, $DIAM_{ha}^{rel}$ drops with decreasing
$\partial_{rat}^{min}(1)$. This is not the case when we decrease
$S(1)$, as shown in table 6.\\

\begin{tabular}{|c|}
\hline Table 6: \textbf{Triadic Snowflake} ($S(1)$ decreases) \label{table6}\\
\hline
\end {tabular}

\vspace{.1in}

\begin{tabular}{|c|c|c|c|c|c|c|}
\hline $\partial_{rel}(1)$ & $S_{rel}(1)$ & $PR(\infty)$ &
$DIAM(\infty)$ & $D_{ha}^{rel}(\infty)$ & $S_{ha}^{rel}(\infty)$ &
$DIAM_{ha}^{rel}(\infty)$ \\
\hline $6.25*10^{-4}$ & 0.4277 & 214.5546 & 1.112 & 0.0348 &
0.5263
& 0.1964\\
\hline 0.3973 & 0.231 & 169.2816 & 0.7011 & 0.0144 & 0.4883 &
0.5033\\
\hline 0.4776 & 0.063 & 156.8692 & 0.5492 & 0.0153 & 0.4897 &
0.7333\\
\hline 0.4752 & 0.063 & 319.6411 & 0.5881 & 0.0163 & 0.4897 &
0.7852\\
\hline 0.4895 & 0.05 & 222.6769 & 0.5999 & 0.025 & 0.5310
&
0.6707\\
\hline 0.4969 & 0.042 & 222.6763 & 0.6102 & 0.025 & 0.531 & 0.6707\\
\hline
\end{tabular} \\

\noindent Unlike in the case of the square Snowflake, we notice
here that a decrease in the initial node entropy results in a
decrease in $DIAM(\infty)$.

\newpage
\noindent We next present in tables 7 and 8 some random walks
started at high entropy/high $\partial_{rat}^{min}$ nodes on both
domains. We point out that on both domains, no node with
$\partial_{rat}^{min}(1)$ and $S(1)$ greater than $1/2$ of their
respective maximum values exists.

\begin{tabular}{|c|}
\hline Table 7: \textbf{Square Snowflake} \label{table7}\\ \hline
\end {tabular}

\vspace{.1in}

\begin{tabular}{|c|c|c|c|c|c|c|}
\hline $\partial_{rel}(1)$ & $S_{rel}(1)$ & $PR(\infty)$ &
$DIAM(\infty)$ & $D_{ha}^{rel}(\infty)$ & $S_{ha}^{rel}(\infty)$ &
$DIAM_{ha}^{rel}(\infty)$ \\
\hline 0.7109 & 0.3328 & 25.6465 & 1.2798 & 0.0644 & 0.707 &
0.2032\\ \hline 0.7108 & 0.3312 & 25.5997 & 1.2847 & 0.065 &
0.6972 &
0.1961\\
\hline 0.5625 & 0.3394 & 25.762 & 0.9719 & 0.0548 & 0.7071 &
0.5256\\ \hline 0.5624 & 0.3375 & 25.7173 & 0.9733 & 0.0563 &
0.6972 &
0.5315\\
\hline 0.3332 & 0.5414 & 26.6456 & 1.2564 & 0.0538 & 0.7071 &
0.2661\\
\hline 0.3125 & 0.4227 & 25.5996 & 1.0088 & 0.0796 & 0.6972 &
0.4651\\ \hline
\end{tabular} \\

\noindent Entropy values of high amplitude points are the same as
for starting nodes with low $\partial_{rat}^{min}(1)$ and $S(1)$
values. $DIAM_{ha}^{rel}(\infty)$ varies a lot.

\begin{tabular}{|c|}
\hline Table 8: \textbf{Triadic Snowflake}\label{table8} \\ \hline

\end {tabular}

\vspace{.1in}

\begin{tabular}{|c|c|c|c|c|c|c|}
\hline $\partial_{rel}(1)$ & $S_{rel}(1)$ & $PR(\infty)$ &
$DIAM(\infty)$ & $D_{ha}^{rel}(\infty)$ & $S_{ha}^{rel}(\infty)$ &
$DIAM_{ha}^{rel}(\infty)$ \\
\hline 0.575 & 0.0113 & 211.4606 & 0.8316 & 0.0229 & 0.5256 &
0.275\\ \hline 0.527 & 0.3067 & 222.0036 & 0.7875 & 0.0236 &
0.5265 &
0.3438\\
\hline 0.4994 & 0.4138 & 222.0426 & 0.8087 & 0.0236 & 0.5265 &
0.3155\\
 \hline 0.4598 & 0.5141 & 222.0645 & 0.8299 & 0.0236 &
0.5265 &
0.2896\\
 \hline 0.3065 & 0.743 & 86.2551 & 0.9098 & 0.0151 & 0.4939
&
0.1511\\
\hline 0.3248 & 0.5546 & 87.0524 & 0.8781 & 0.026 & 0.4896 &
0.1892 \\
\hline
\end{tabular} \\

\noindent  $DIAM_{ha}^{rel}(\infty)$ values are lower than on the
square domain. $DIAM(\infty)$ also seems to be lower and
$D_{ha}^{rel}(\infty)$ only lies at about $2\%$ of the maximum
value instead of at about $5.5\%$ on the other domain. Again,
participation ratios vary on the triadic domain while they are
constant for the square Snowflake.

\newpage
\noindent Tables 9 and 10 refer to the simulations run with high
$\partial_{rat}^{min}(1)$ and low $S(1)$ values.

\begin{tabular}{|c|}
\hline Table 9: \textbf{Square Snowflake} \label{table9}\\ \hline
\end {tabular}

\vspace{.1in}

\begin{tabular}{|c|c|c|c|c|c|c|}
\hline $\partial_{rel}(1)$ & $S_{rel}(1)$ & $PR(\infty)$ &
$DIAM(\infty)$ & $D_{ha}^{rel}(\infty)$ & $S_{ha}^{rel}(\infty)$ &
$DIAM_{ha}^{rel}(\infty)$ \\
\hline 0.7656 & $1.31*10^{-8}$ & 25.5997 & 1.2649 & 0.0683 &
0.6972 & 0.2045\\
\hline 0.9376 & $6.49*10^{-5}$ & 25.6468 & 0.8764 & 0.0548 &
0.7071 & 0.6714\\
\hline 0.9376 & $2.96*10^{-4}$ & 25.6473 & 0.9209 & 0.0592 &
0.7071 & 0.596\\
\hline 0.9376 & $2.1*10^{-9}$ & 25.6053 & 0.8643 & 0.0501 & 0.6972& 0.702\\
\hline  0.2188 & $4.07*10^{-12}$ & 25.6595 & 0.925 & 0.0533 &
0.7071 & 0.5922\\
\hline 0.9376 & $6.49*10^{-5}$ & 25.5997 & 0.8796 & 0.0563 &
0.6972
& 0.6725\\
\hline
\end{tabular} \\

\begin{tabular}{|c|}
\hline Table 10: \textbf{Triadic Snowflake}\label{table10} \\
\hline

\end {tabular}

\vspace{.1in}

\begin{tabular}{|c|c|c|c|c|c|c|}
\hline $\partial_{rel}(1)$ & $S_{rel}(1)$ & $PR(\infty)$ &
$DIAM(\infty)$ & $D_{ha}^{rel}(\infty)$ & $S_{ha}^{rel}(\infty)$ &
$DIAM_{ha}^{rel}(\infty)$ \\
\hline 0.4952 & 0.04 & 222.6771 & 0.6844 & 0.025 & 0.531 &
0.4638\\
\hline 0.59 & $4.32*10^{-3}$ & 221.6813 & 0.6529 & 0.0182 &
0.5228 & 0.5147\\
\hline 0.6383 & 0.6554 & 222.6791 & 0.7485 & 0.0253 & 0.5324 & 0.3408\\
\hline 0.734 & 0.147 & 210.968 & 0.747 & 0.0219 & 0.524 & 0.3908
\\
\hline 0.772 & 0.1176 & 213.7942 & 0.7185 & 0.0224 & 0.5255 &
0.4135 \\
\hline 0.8483 & $1.05*10^{-6}$ & 211.4231 & 0.717 & 0.0229 &
0.5256
& 0.3808\\
\hline
\end{tabular} \\

\noindent Note that low $S(1)$ values are not sufficient for
$DIAM_{ha}^{rel}(\infty)$ values to drop.

\newpage
\noindent Finally, tables 11 and 12 provide the results of the
simulations run with low $\partial_{rat}^{min}(1)$ and high $S(1)$
values.

\begin{tabular}{|c|}
\hline Table 11: \textbf{Square Snowflake}\label{table11} \\
\hline

\end {tabular}

\vspace{.1in}

\begin{tabular}{|c|c|c|c|c|c|c|}
\hline $\partial_{rel}(1)$ & $S_{rel}(1)$ & $PR(\infty)$ &
$DIAM(\infty)$ & $D_{ha}^{rel}(\infty)$ & $S_{ha}^{rel}(\infty)$ &
$DIAM_{ha}^{rel}(\infty)$ \\
\hline 0.1876 & 0.9554 & 25.5997 & 1.3179 & 0.0624 & 0.697 &
0.1221\\
\hline 0.1388 & 0.949 & 25.6047 & 0.9904 & 0.0563 & 0.6972 &
0.5015\\
\hline 0.0936 & 0.9363 & 25.5997 & 1.3681 & 0.0692 & 0.6972 &
0.087\\
\hline 0.09 & 0.9358 & 25.6455 & 1.3514 & 0.0548 & 0.7071 &
0.0842\\
 \hline 0.078 & 0.9044 & 25.6456 & 1.3486 & 0.0548 & 0.7071
&
0.0882\\
\hline 0.078 & 0.9044 & 25.5996 & 1.3208 & 0.0579 & 0.6972 & 0.1165\\
\hline
\end{tabular} \\

\noindent Note the high and low values of the relative diameter of
high amplitude points of row 2 an 3, respectively.\\

\begin{tabular}{|c|}
\hline Table 12: \textbf{Triadic Snowflake} \label{table12}\\
\hline
\end {tabular}

\vspace{.1in}

\begin{tabular}{|c|c|c|c|c|c|c|}
\hline $\partial_{rel}(1)$ & $S_{rel}(1)$ & $PR(\infty)$ &
$DIAM(\infty)$ & $D_{ha}^{rel}(\infty)$ & $S_{ha}^{rel}(\infty)$ &
$DIAM_{ha}^{rel}(\infty)$ \\
\hline 0.206 & 0.7478 & 86.5501 & 0.8474 & 0.0292 & 0.4897 &
0.2385\\
\hline 0.1618 & 0.9957 & 86.6407 & 0.8497 & 0.0313 & 0.4897 &
0.2546\\
\hline 0.1147 & 0.8944 & 222.1414 & 0.8283 & 0.0238 & 0.5282 &
0.3386\\
\hline 0.09 & 0.8361 & 213.3961 & 0.8561 & 0.0236 & 0.5254 &
0.1841\\
\hline 0.082 & 0.7879 & 86.8069 & 0.8565 & 0.0312 & 0.4897 &
0.275\\
\hline 0.061 & 0.7016 & 213.8787 & 0.9734 & 0.0233 & 0.5243 &
0.121\\
 \hline
\end{tabular} \\

\noindent Even though we start in a high entropy region, it seems
to be sufficient to have low $\partial_{rel}(1)$ for a drop in
$DIAM_{ha}^{rel}(\infty)$ to occur.

\section{Results and Discussion}

\noindent Our numerical calculations show the following results on
both domains:

\begin{itemize}
\item All five quantities of interest approach an asymptotic
limit. \item $D_{ha}^{rel}(\infty)$ and $S_{ha}^{rel}(\infty)$ are
independent of the starting node. Values lie at about $5.5\%$ and
$70\%$ of the respective maximum values on the square Snowflake
and at about $2.5\%$ and $50\%$ of the respective maximum values
on the triadic Snowflake.\item The asymptotic limit of the
Diameter ($DIAM(\infty)$) varies a lot within a region but overall
seems to be a little lower when walks are started in low
distance/low entropy regions. \item When $S(1)$ and
$\partial_{rat}^{min}(1)$ are very low, the diameter of high
amplitude points ($DIAM_{ha}^{rel}(\infty)$) is significantly
smaller than for walks started in any other region. \item No node
with $S(1)$ and $\partial_{rat}^{min}(1)$ values larger than half
their respective maximum values can be found.
\end{itemize}

\noindent Comparing results on the triadic and square domains, we
observe the following:

\begin{itemize}
\item Astonishingly, the participation ratio of the square
Snowflake is approximately the same for any starting node. On the
triadic domain, we find three different approximate values: even
if $S(1)$ and $\partial_{rat}^{min}(1)$ values are similar (row
three and four in table \ref{table6}), we can observe different
participation ratios. \item The participation ratio on the triadic
Snowflake is higher than on the square Snowflake: On the square
domain the participation ratio lies at $0.5\%$, the three values
on the triadic domain correspond to about 1, 3 and $4\%$ of the
total number of nodes in the discretization. \item On the triadic
domain, no node with
$\partial_{rat}^{min}<\frac{1}{2}\max{(\partial_{rat}^{min})}$ and
$S<\frac{1}{24}\max{(S)}$ can be found, while on the square
domain, we can find a node satisfying $S<\frac{1}{251}\max{(S)}$
for
$\partial_{rat}^{min}<\frac{1}{2}\max{(\partial_{rat}^{min})}$.
\item The minimum ratio distances of high entropy points on the
square domain are constant and lie at about $5\%$ of its maximum
value on the domain. On the triadic Snowflake, minimum ratio
distances of high entropy points vary: we find values from $1\%$
to $3\%$ of the maximum. These seem to be independent of where the
starting node lies. \item For both domains, the
$S_{ha}^{rel}(\infty)$ values do not show sensitivity to initial
conditions. It is rather surprising, though, that they are larger
on the square Snowflake.
\end{itemize}

\noindent In summary, we can say that random walks started at low
$\partial_{rat}^{min}$ nodes show localization in the form of
decreased diameter of high amplitude points. This localization
occurs on both domains with the same frequency. The value of the
initial node entropy does not seem to have a great influence on
the dynamics of the system: a decrease does not necessarily result
in a decrease of the diameter of high amplitude points.
$D_{ha}^{rel}(\infty)$ and $S_{ha}^{rel}(\infty)$ always level off
to about the same value on both domains, with the
$D_{ha}^{rel}(\infty)$ very low and  $S_{ha}^{rel}(\infty)$ values
relatively large. This is a little surprising and requires further
investigation. The main result of the numerical simulations is
that the dynamics of the two systems are very similar. Hence, we
can not conclude that narrow channels in the geometry of the
square domain cause localization of the Dirichlet eigenfunctions.
But there certainly are noticeable geometric differences on the
two domains. The nodal entropies found on the square Snowflake are
significantly higher. Maybe the greater diversity in
$\partial_{rat}^{min}(1)$ distribution on this domain is the cause
for localization found in some of the Dirichlet eigenfunctions.
This hypothesis should be further investigated. One should create
weighted graphs with different nodal entropy distributions and
test how the dynamics of the system is influenced. Another
direction of further research in this area would be to conduct a
thorough analysis of the transfer matrix $\mathbf{T}$. It would be
interesting to determine if the nature of the eigenvalues and
eigenfunctions of $\mathbf{T}$ reveals information about the
course taken by a random walker.

\begin{center}
\textbf{Acknowledgements}\\
\noindent The authors would like to thank Erin Pearse for fruitful
conversations.
\end{center}

\end{document}